\documentclass[12pt]{amsart}
\usepackage[pdftex]{graphicx}
\usepackage{amssymb}
\usepackage{MnSymbol}
\usepackage{mathtools}

\usepackage[mathscr]{eucal}
\usepackage{amssymb, amsmath,array, amscd}
\usepackage{enumerate}
\usepackage{graphicx}
\usepackage{url}
\usepackage[colorlinks,plainpages,backref]{hyperref}

\usepackage{tikz}
\usepackage{pifont}
\usetikzlibrary{matrix}
\usepackage{amssymb}
\usepackage{MnSymbol}
\usepackage{mathtools}

\newtheorem{thm}{Theorem}[section]
\newtheorem{defn}{Definition}[section]
\newtheorem{lem}{Lemma}[section]
\newtheorem{cor}{Corollary}[section]
\newtheorem{ex}{Example}[section]
\newtheorem{rem}{\it Remark}[section]

\begin{document}

\title{Graphs of Schemes Associated to  Group Actions}

\author{Al\.i Ula\c{s} \"Ozg\"ur K\.i\c{s}\.isel}
\address{Middle East Technical University \\06531 Ankara \\  Turkey} 
\email{akisisel@metu.edu.tr}

\author{Eng\.in \"Ozkan}
\address{S\"uleyman Demirel University \\ Isparta \\ Turkey}
\email{enginozkan@sdu.edu.tr}


\begin{abstract}
Let $X$ be a proper algebraic scheme over an algebraically closed field. We assume that a torus $T$ acts on $X$ such that the action has isolated fixed points. The $T$-graph of $X$ can be defined using the fixed points and the one dimensional orbits of the $T$-action. If the upper Borel subgroup of the general linear group with maximal torus $T$ acts on $X$, then we can define a second graph associated to $X$, called the $A$-graph of $X$. We prove that the $A$-graph of $X$ is connected if and only if $X$ is connected. We use this result to give a proof of Hartshorne's theorem on the connectedness of Hilbert scheme in the case of $d$ points in $\mathbb{P}^{n}$. 
\end{abstract}

\maketitle

\section{Introduction}

Let $X$ be a proper algebraic scheme over an algebraically closed field $k$ and let $T=(G_m)^{n}=(k^{*})^{n}$ be an $n$-dimensional algebraic torus. In the case that $X$ admits a $T$-action with isolated fixed points one can extract substantial information about $X$ by looking at these fixed points, through equivariant localization or other means (see for instance \cite{BBCM}, \cite{GKM}, \cite{Ive}). In \cite{AS} the concept of a $T$-graph is defined for the Hilbert scheme and its multigraded versions. This definition can be generalized to any setting where the $T$-action has isolated  fixed points. In \cite{HM} further properties of $T$-graphs are obtained and some explicit computations are given. 

If $X$ admits an action of a larger group, for instance a Borel subgroup of the general linear group such that the induced maximal torus action has isolated fixed points, then the information that one can extract about $X$ often becomes substantially stronger. As an example, the results of \cite{AC} imply that if $X$ admits a compatible $(G_{m},G_{a})$-action such that the $G_{a}$-fixed locus consists of a single point, then the Poincar\'{e} polynomial of $X$ is a product of cyclotomic polynomials. 

In this paper, we first generalize the definition of a $T$-graph to an arbitrary proper algebraic scheme $X$ admitting a $T$-action with isolated fixed points and show that $X$ is connected if and only if its $T$-graph is connected. These issues are discussed in section 2. In section 3, assuming that $X$ admits a Borel group action, we define the notion of an $A$-graph, which is obtained from the fixed locus of $X$ under the action of a principal nilpotent subgroup. The $A$-graph of $X$ is often substantially smaller compared to the $T$-graph. For instance, in the setting described in the last sentence of the previous paragraph, the $A$-graph consists of just a single vertex. We prove that $X$ is connected if and only if its $A$-graph is connected. 

In section 4, we give a proof of Hartshorne's theorem of the connectedness of the Hilbert scheme of $\mathbb{P}^{n}$ (\cite{H}) in the case where the Hilbert polynomial is a constant $d$, by explicitly computing a spanning subtree of the $A$-graph of the Hilbert scheme. This proof can be generalized to the case of an arbitrary Hilbert polynomial, but this is not given here. The idea of using Borel fixed monomial ideals is present  both in Hartshorne's proof \cite{H} and A. Reeves' proof \cite{Ree}; the current proof provides a somewhat different organization. Computation of the full $A$-graphs of Hilbert schemes of points for small values of  $d$ and $n$ will be the subject of a future paper. 

\section{$T$-Graph of a Scheme}

Let $k$ be an algebraically closed field and $X$ a proper algebraic scheme over $k$. Let $T=(G_m)^{n}=(k^{*})^{n}$ be an $n$-dimensional algebraic torus. Assume that $T$ acts on $X$ such that the fixed point scheme $X^{T}$ consists of isolated points. 

\begin{defn}
The graph whose vertices correspond to elements of $X^{T}$ and having an edge between $p,q\in X^{T}$ if and only if there exists a $1$-dimensional $T$-orbit $l$ whose closure contains both $p$ and $q$ is called the $T$-graph of $X$. 
\end{defn}

The $T$-graph of the Hilbert scheme of $\mathbb{P}^n$ was defined in \cite{AS}. The $T$-graphs of multigraded Hilbert schemes of $\mathbb{P}^{n}$ are also implicitly present in the discussion in  \cite{AS} although the emphasis there is rather on the induced subgraphs obtained from the ordinary grading. This definition is given in terms of Gr\"{o}bner degenerations, but it is equivalent to the one above if we assume that the action on this Hilbert scheme is induced from the standard action of the maximal torus of $GL(n+1)$ on $\mathbb{P}^{n}$. Some details are provided in the two examples below.   

\begin{ex} 
Let $X=\mathbb{P}^{n}$ and $T$ be the maximal torus of $GL(n+1)$ containing diagonal matrices $diag(d_{0},\ldots,d_{n})$ such that $d_{i}\neq 0$. Consider the standard action of $GL(n+1)$ on $X$ and the corresponding action of $T$ on $X$ given by $diag(d_{0},\ldots,d_{n})\cdot [x_{0}:\ldots:x_{n}]= [d_{0}x_{0}:\ldots:d_{n}x_{n}]$. The fixed points of this action are the $n+1$ points of the form $p_{i}=[0:\ldots:0:1:0:\ldots:0]$ for $i=0,1,\ldots,n$ where the only nonzero entry of $p_{i}$ is its $i$th coordinate. The open subset of the line $l_{ij}=[0:\ldots:\lambda: \ldots: \mu:\ldots:0]$, where the only nonzero entries $\lambda, \mu$ are at the $i^{th}$ and $j^{th}$ coordinates respectively, is a $1$-dimensional $T$-orbit. Since $p_i$ and $p_j$ lie in the closure of $l_{ij}$, any two vertices are connected by an edge. Hence the $T$-graph of $X$ for this action is the $1$-skeleton of the standard $n$-simplex. 
\end{ex} 

\begin{ex}
Let $X=Hilb(\mathbb{P}^n)$, the Hilbert scheme of $\mathbb{P}^{n}$. Elements of $X$ correspond to homogenous ideals in $k[x_{0},x_{1},\ldots,x_{n}]$. The $T$-action on $\mathbb{P}^{n}$ induces a $T$-action on $X$. It is easy to see that the fixed points of this action are precisely the monomial ideals. Indeed, contrary to the claim, assume that the ideal $I$ is a fixed point but it is not a monomial ideal. Then there exists a polynomial $P(x_{0},\ldots, x_{n})$ in $I$ so that none of its monomial summands are in $I$. Suppose that $P$ is a sum of $r$ monomials. We can assign integer weights $m_{0},\ldots,m_{n}$ to the variables $x_{0},\ldots,x_{n}$ such that each monomial summand of $P$ has a different weight. Act on $P$ by elements of the form $diag(t^{m_{0}},\ldots,t^{m_{n}})$ for $r$ different values of $t$. The resulting polynomials $P_1,\ldots,P_r$ must be in $I$ because $I$ is a fixed point of the action. Then an elementary computation involving a Van der Monde determinant shows that the set $\{P_1,\ldots,P_r\}$ is linearly independent. This implies that each monomial in $P$ can be written as a linear combination of $\{P_{1},\ldots,P_{r}\}$, contradicting the initial assumption. This contradiction shows that the only fixed points of the $T$-action are monomial ideals. Thus we determined the vertices of the $T$-graph. Determining the edges of the $T$-graph is a more difficult task, see  \cite{AS} or \cite{HM} for a discussion.   
\end{ex} 

\begin{thm} \label{thm:Tgraph} Suppose that $T$ acts on the proper algebraic scheme $X$ with isolated fixed points. Then $X$ is connected if and only if its $T$-graph is connected. 
\end{thm}

\begin{rem} The analogous statement is false in the category of smooth manifolds. For example, take $k=\mathbb{C}$ and assume that $X$ is $S^{2}$, the $2$-sphere. Identify $S^{2}$ with the Riemann sphere so that the North pole corresponds to $\infty$, the South pole corresponds to $0$ and $1$ is on the equator. Let $S^{1}\times S^{1}$ act on $X$ in the following way: The first factor acts by rotations parallel to the equator. For the second factor, choose a M\"{o}bius transformation $M_{\alpha}$ sending $0$ to $0$, $\infty$ to $\infty$ and $1$ to $\alpha$ where $\alpha$ is not on the equator. Then assume that the action of the second factor is the conjugation of the first action by $M_{\alpha}$. Then the fixed points of both actions are the north and south poles. The $1$-dimensional orbits of the two actions are two families of circles. A circle in one family intersects circles of the other family transversely on a dense subset. This implies that the $S^{1}\times S^{1}$-action on $X$ has two fixed points, no $1$-dimensional orbits and a single $2$-dimensional orbit. We can turn this into a $T$-action by retracting $T=\mathbb{C}^{*}\times \mathbb{C}^{*}$ radially to $S^{1}\times S^{1}$ and composing the retraction with the action described above. The $T$-graph of $X$ consists in just two disconnected vertices and no edges. Hence the statement is not true in this case. 
\end{rem}

\begin{lem} \label{lem:parameter}

(i) There exists a $1$-parameter subgroup $G_m\hookrightarrow T$such that $X^{T}=X^{G_m}$. 

(ii) There exists a $2$-parameter subgroup $G_{m}\times G_{m}\hookrightarrow T$ such that the $T$-graph of $X$ agrees with the $G_m\times G_m$-graph of $X$. 
\end{lem} 

\noindent\textit{Proof} It is enough to prove both statements in the case where $X$ is irreducible. 

(i) By the Bialynicki-Birula decomposition \cite{BB1} there exists a big cell $U\subset X$ that contains a unique fixed point $a$ and the $T$-action restricted to $U$ is conjugate to a linear action of $T$ on a vector space $V$. Since the action has a unique fixed point, all weights of this action must be nontrivial. Suppose that we choose a basis that diagonalizes the action, namely $t=(t_{1},\ldots,t_{n})\in T$ acts on $V$ by $diag(\lambda_{1}(t),\ldots, \lambda_{N}(t))$. Assume that $\lambda_{i}(t)=t_{1}^{e_{i1}}t_{2}^{e_{i2}}\ldots t_{n}^{e_{in}}$. Since there are no trivial weights, for each $i$ there exists $j$ such that $e_{ij}\neq 0$. Now consider an arbitrary one-parameter subgroup $\varphi: G_{m}\hookrightarrow T$ so that $\varphi(s)=(s^{a_{1}},s^{a_{2}},\ldots,s^{a_{n}})$. The weights restricted to this subgroup are nontrivial if and only if for each $i$, the sum $\sum_{j=1}^{n}e_{ij}a_{j}$ is nonzero. Since this is a set of finitely many nontrivial inequalities in integers, there exists $(a_{1},\ldots,a_{n})$ 
such that $U^{T}=U^{G_{m}}$. Now replace $X$ by $X-U$ and proceed inductively. Since at every step finitely many equalities are to be avoided and there are finitely many steps, there exist one-parameter subgroups $\varphi:G_{m}\hookrightarrow T$ such that $X^{G_{m}}=X^{T}$. 

(ii) Let $U$ and $V$ be as in part (i) of the proof. The $T$-graph of $X$ has a unique vertex $v$ on $U$ corresponding to the unique fixed point. The edges of the graph emanating from $v$ correspond to lines spanned by eigenvectors of the $T$-action on $V$. Choose a basis as in part (i). A two parameter subgroup $\varphi:G_{m}\times G_{m}\hookrightarrow T$ induces the same edges as $T$ if and only if no new vertices are produced and the eigenspaces of the $G_{m}\times G_{m}$-action are the same as the eigenspaces of the $T$-action. Say $\varphi(s,u)=(s^{a_{1}}u^{b_{1}},s^{a_{2}}u^{b_{2}},\ldots,s^{a_{n}}u^{b_{n}})$. Then $\lambda_{i}(s,u)=s^{\sum_{j=1}^{n} a_{j}e_{ij}}u^{\sum_{j=1}^{n} b_{j}e_{ij}}$. Therefore the $T$-graph restricted to $U$ agrees with the $G_{m}\times G_{m}$-graph restricted to $U$ if and only if the following inequalities are satisfied for all pairs $i,k$ with $(e_{i1},\ldots,e_{in})\neq (e_{k1},\ldots,e_{kn})$:

\begin{equation*}
(\sum_{j=1}^{n} a_{j}e_{ij}, \sum_{j=1}^{n} b_{j}e_{ij})\neq (\sum_{j=1}^{n} a_{j}e_{kj}, \sum_{j=1}^{n} b_{j}e_{kj}) 
\end{equation*}
and $(\sum_{j=1}^{n} a_{j}e_{ij}, \sum_{j=1}^{n} b_{j}e_{ij})\neq (0,0)$ for all $i$.
Since a finite set of nontrivial equalities need to be avoided, there exist solutions for $a_{i},b_{i}$. The rest of the argument proceeds  as in (i), by considering $X-U$ and doing induction. $\Box$    

We will adapt the strategy in \cite{BB2} to our case in order to prove theorem \ref{thm:Tgraph}. We now recall a definition from \cite{BB2} : 
\begin{defn}
Let $A, X_{1}, X_{2}$ be schemes. We say that $X_{1}$ is simply $A$-equivalent to $X_{2}$ if there exists an isomorphism $X_{1}-X_{2}\rightarrow A\times Y$ for some scheme $Y$. The equivalence closure of simple $A$-equivalence is called $A$-equivalence. 
\end{defn} 

Let $K$ be a finite field so that $char(K)\neq char(k)$. Let $H^{i}(X)$ denote the $i$-th \'{e}tale cohomology group of $X$ with proper supports, with coefficients in the constant sheaf $K$. 

\begin{lem} \label{lem:cohomology}
Say $X_{1}$ and $X_{2}$ are schemes such that $X_{1}$ is $G_{m}\times G_{m}$-equivalent to $X_{2}$. Then $H^{0}(X_{1})\cong H^{0}(X_{2})$.
\end{lem}

\noindent \textit{Proof} Since $H^{0}(G_{m})=0$ and $H^{1}(G_{m})=K$, by the K\"{u}nneth formula we obtain  $H^{0}(G_{m}\times G_{m})=H^{1}(G_{m}\times G_{m})=0$. 

It is enough to prove the statement when $X_{1}$ is simply $G_{m}\times G_{m}$-equivalent to $X_{2}$. Hence we may assume that $X_{1}-X_{2}\cong U = G_{m}\times G_{m} \times U_{1}$. If $U_{1}$ is singular then we can stratify its singular locus and break this equivalence into finitely many steps in each of which $U_{1}$ is nonsingular. Hence we may assume that $U_{1}$ is nonsingular. By the K\"{u}nneth formula, there exists a spectral sequence such that $E_{2}^{p,q}=H^{p}(U_{1})\otimes H^{q}(G_{m}\times G_{m})$ converging to $H^{p+q}(U)$. 

First, let us consider the case of $dim(U)=2$. Then $U$ is isomorphic to a disjoint union of copies of $G_{m}\times G_{m}$. It suffices to consider one copy at a time, so assume that $U\cong G_{m}\times G_{m}$. If $\overline{U}=U$ then $H^{0}(\overline{U})=H^{0}(U)=0$. The excision sequence in cohomology says that $0\rightarrow H^{0}(X_{1}-\overline{U})\rightarrow H^{0}(X_{1})\rightarrow H^{0}(\overline{U})$ is exact. But $X_{1}-\overline{U}=X_{2}$ and $H^{0}(\overline{U})=0$. Therefore $H^{0}(X_{1})\cong H^{0}(X_{2})$. If $\overline{U}\neq U$ then first consider the excision sequence with respect to the inclusion $U\hookrightarrow \overline{U}$. Then 

\[ H^{0}(U)\rightarrow H^{0}(\overline{U})\rightarrow H^{0}(\overline{U}-U)\rightarrow H^{1}(U) \]
is exact. But $H^{0}(U)=H^{1}(U)=0$, therefore $H^{0}(\overline{U})\cong H^{0}(\overline{U}-U)$. Now, consider the excision sequences with respect to the inclusions $\overline{U}\hookrightarrow X_{1}$ and $\overline{U}-U \hookrightarrow X_{1}-U$. 

\hspace{-1cm}\begin{tikzpicture}
  \matrix (m) [matrix of math nodes,row sep=3em,column sep=3em,minimum width=2em]
  {
     0 & H^{0}(X_{1}-\overline{U}) &H^{0}(X_{1}) & H^{0}(\overline{U}) & H^{1}(X_{1}-\overline{U}) \\
     0 & H^{0}(X_{1}-\overline{U}) & H^{0}(X_{1}-U) & H^{0}(\overline{U}-U) & H^{1}(X_{1}-\overline{U}) \\};
  \path[-stealth]
    (m-1-1) edge (m-1-2)
    (m-1-2) edge (m-1-3)  
    (m-1-3) edge (m-1-4) 
     edge node [right] {$\alpha$} (m-2-3)
    (m-1-4) edge (m-1-5) 
    edge node [right] {$\cong$} (m-2-4)
    (m-2-1.east|-m-2-2) edge (m-2-2)
    (m-2-2) edge (m-2-3)
    (m-2-3) edge (m-2-4) 
    (m-2-4) edge (m-2-5) 
    (m-1-5) edge node [right] {$\cong$} (m-2-5)      
    (m-1-2) edge node [right] {$\cong$} (m-2-2)           ;
\end{tikzpicture}
By the $5$-lemma, the map $\alpha: H^{0}(X_{1})\rightarrow H^{0}(X_{1}-U)=H^{0}(X_{2})$ is an isomorphism, which finishes the proof in this case. 

Next, consider the case of $dim(U)>2$. We may assume that none of the connected components of $U_{1}$ is complete. Else, stratify $U_{1}$ and subdivide the equivalence into more steps, where the stratum with $dim(U_{1})=0$ is taken care of as above. 
Hence $H^{0}(U_{1})=0$. Since $E_{2}^{p,q}(U)=H^{p}(U_{1})\otimes H^{q}(G_{m}\times G_{m})$ we get $E_{2}^{0,0}(U)=E_{2}^{0,1}(U)=E_{2}^{1,0}(U)=0$. This implies $H^{0}(U)=H^{1}(U)=0$. The rest of the argument is identical to the one in the first case above. $\Box$

\noindent \textit{\textbf{Proof of Theorem \ref{thm:Tgraph}}} Let $Y$ denote the union of $0$ and $1$-dimensional orbits of the $T$-action on $X$. By lemma \ref{lem:parameter} part (ii), there exists a subgroup of $T$ isomorphic to $G_{m}\times G_{m}$ whose $0$ and $1$-dimensional orbits agree with those of $T$. It is clear that $X$ is $G_{m}\times G_{m}$-equivalent to $Y$. Hence, by lemma \ref{lem:cohomology}, $H^{0}(X)\cong H^{0}(Y)$. Since $X$ is proper, $X$ is connected if and only if $H^{0}(X)\cong K$ and the same holds for $Y$. On the other hand, by definition, $Y$ is connected if and only if the $T$-graph of $X$ is connected. This establishes the claim. $\Box$  

\section{Additive Group Actions and $A$-Graphs} 

Let $G_{a}$ denote the 1-dimensional connected
additive algebraic group. 

 \begin{defn} Let $X$ be an algebraic scheme over $k$ having $G_{m}$ and $G_{a}$-actions such that
\begin{eqnarray*} 
\lambda: G_{m}\times X \rightarrow X, \qquad ((t,x)\mapsto
\lambda(t)\cdot x) \\
\theta: G_{a}\times X \rightarrow X, \qquad ((u,x)\mapsto
\theta(u)\cdot x)
\end{eqnarray*}
The actions $\lambda$ and $\theta$ are said to be compatible if there exists an integer $p\geq 1$ such that
\[\lambda(t)\cdot \theta(u)\cdot \lambda(t^{-1})=\theta(t^p\cdot u)\]
for all $t\in G_{m}$ and $u\in G_{a}$.
\end{defn}

Suppose that $X$ admits compatible $G_m$ and $G_{a}$-actions. Then the $G_{a}$-fixed locus $X^{G_{a}}$ of $X$ is $G_{m}$-invariant. Indeed, suppose that $x\in X^{G_{a}}$. Then $\theta(u) \lambda(t)x=\lambda(t)\theta(t^{-p}u)x=\lambda(t)x$. Thus $\lambda(t)x\in X^{G_{a}}$. 

Let $E_{i,j}$ denote the elementary matrix whose $ij$-entry is $1$ and all others $0$. Say $E=\sum_{i=1}^{n-1} E_{i,i+1}$. Then a copy of $G_{a}$ consisting in matrices of the form $exp(aE)$ for $a\in k$ sits inside $GL(n)$. Let $T$ be the the maximal torus of $GL(n)$ consisting of diagonal matrices. Then it is easy to check that $1$-parameter subgroups $G_{m}$ of $T$ compatible with this $G_{a}$ are subgroups containing matrices of the form: 
\[ diag(t^c, t^{c-p}, t^{c-2p},\ldots, t^{c-(n-1)p}).  \]

Say $X$ is a proper algebraic scheme over $k$. Let $B$ be  the upper Borel subgroup of $GL(n)$ and $T$ the maximal torus in $B$ consisting of diagonal matrices. Suppose that $B$ acts on $X$ such that the induced action of $T$ on $X$ has isolated fixed points. Let $G_{a}$ be the subgroup of $B$ consisting in the matrices $exp(aE)$ as above. The following lemma is a refinement of lemma \ref{lem:parameter}: 

\begin{lem}\label{lem:compat}
There exists a $2$-parameter subgroup $G_{m}\times G_{m}\hookrightarrow T$ such that 

(i) the $T$-graph of $X$ agrees with the $G_{m}\times G_{m}$-graph of $X$,

(ii) actions of both $G_{m}$ factors are compatible with the $G_{a}$-action. 
\end{lem}

\noindent \textit{Proof} In order for condition (ii) to be satisfied, the two $G_m$ factors should contain matrices of the form $diag(t^{a_{1}},\ldots,t^{a_{n}})$ and $diag(t^{b_{1}},\ldots,t^{b_{n}})$ respectively where $a_{i}=c_{1}-(i-1)p_{1}$ and $b_{i}=c_{2}-(i-1)p_{2}$ for some $c_{1}, c_{2}, p_{1}, p_{2}$. One can check that the inequalities $(\sum_{j=1}^{n} a_{j}e_{ij}, \sum_{j=1}^{n} b_{j}e_{ij})\neq (\sum_{j=1}^{n} a_{j}e_{kj}, \sum_{j=1}^{n} b_{j}e_{kj})$ and $(\sum_{j=1}^{n} a_{j}e_{ij}, \sum_{j=1}^{n} b_{j}e_{ij})\neq (0,0)$ in the proof of lemma \ref{lem:parameter} are still nontrivial. Therefore there exist choices of $c_{1},c_{2},p_{1},p_{2}$ such that all inequalities are satisfied. $\Box$ 

\begin{defn} 
Suppose that $X$ is a proper scheme admitting a Borel group action and let $G_{m}\times G_{m}$ and $G_{a}$ be as in \ref{lem:compat}. Then 
the $G_{m}\times G_m$ graph of $X^{G_{a}}$ is called the $A$-graph of $X$ associated to this data. 
\end{defn} 

Recall Horrocks' theorem \cite{H} saying that the fixed locus $X^{S}$ of the action of a solvable group $S$ on a  proper scheme $X$ is connected if and only if $X$ itself is connected. This holds in particular for $G_{a}$-actions. 

\begin{lem} \label{lem:independence}
The $A$-graph of $X$ is independent of the choice of the $G_{m}\times G_{m}$ in its definition. Its vertices are the vertices of the $T$-graph which are $G_{a}$-fixed and its edges are the edges of the $T$-graph which are $G_{a}$-invariant. 
\end{lem} 

\noindent \textit{Proof} The vertices of the $A$-graph are $G_{m}\times G_{m}$-fixed points in $X^{G_{a}}$. By the choice of $G_{m}\times G_{m}$, these are precisely the $T$-fixed points in $X$ which are also $G_{a}$-fixed. Similarly, $1$-dimensional $T$-orbits agree with $1$-dimensional $G_{m}\times G_{m}$-orbits. Say $l$ is such an orbit whose closure contains the fixed points $P$ and $Q$. Suppose that $l$ is $G_{a}$-invariant. But in this case the closure of $l$, which is isomorphic to $\mathbb{P}^{1}$ admits a $G_{a}$-action with at least two fixed points $P$ and $Q$. By Horrocks' theorem, the $G_{a}$-fixed locus of this action on $\mathbb{P}^{1}$ must be connected, therefore all of $l$ must lie in the $G_{a}$-fixed locus $X^{G_{a}}$. Therefore the edges of the $A$-graph are precisely the $1$-dimensional $T$-orbits which are $G_{a}$-invariant or equivalently $G_{a}$-fixed. $\Box$ 

\begin{thm}\label{thm:Agraph}
Suppose that the Borel subgroup $B$ of $GL(n)$ acts on the proper algebraic scheme $X$ such that the induced $T$-action has isolated fixed points. Then the $A$-graph of $X$ is connected if and only if $X$ is connected. 
\end{thm}

\noindent\textit{Proof} By Horrocks' theorem, $X$ is connected if and only if $X^{G_{a}}$ is connected. Applying theorem \ref{thm:Tgraph} 
to $X^{G_{a}}$ we see that $X^{G_{a}}$ is connected if and only if the $A$-graph of $X$ is connected. The result follows. $\Box$ 

\section{$\mathbf{A}$-Graph of $Hilb^{d}(\mathbb{P}^n)$}

In this section we will show that the $A$-graph of $Hilb^{d}(\mathbb{P}^n)$ is connected for any $d>0$ and $n>0$. Given theorem \ref{thm:Agraph}, this will provide us with a proof of Hartshorne's theorem on the connectedness of the Hilbert scheme \cite{RH} in the restricted case of a constant Hilbert polynomial.  Instead of determining the $A$-graph completely, we will explicitly describe the vertices and a spanning subtree of this graph. The case of $n=2$ was treated in \cite{E}.  The arguments here can be generalized to any Hilbert polynomial $p(t)$, but we do not carry this out here.  

Consider the standard representation of $GL(n+1)$ on $K^{n+1}$ by matrix multiplication. This induces a $GL(n+1)$-action on $\mathbb{P}^{n}(K)$ and consequently the upper Borel subgroup $B$ of $GL(n+1)$ acts on $\mathbb{P}^{n}$. Selecting the subgroup $G_{a}$ of $B$ as matrices of the form $exp(aE)$ as in the previous section, the induced $G_{a}$-action $\theta$ on $\mathbb{P}^{n}$ is given by 
\[ \theta(a)([X_{0}:\ldots:X_{n}])=[X_{0}+aX_{1}+\ldots+\frac{a^n}{n!}X_{n}: X_{1}+aX_{2}+\ldots+\frac{a^{n-1}}{(n-1)!}X_{n}:\ldots:X_{n}].  \]
In particular, the $G_{a}$-action on $\mathbb{P}^{n}$ has a unique fixed point $[1:0:\ldots:0]$. 
   
Consider the induced $B$-action on $K[X_{0},\ldots,X_{n}]$. It takes homogenous ideals to homogenous ideals and the action is flat. Therefore, there is an induced $B$-action on the Hilbert scheme $Hilb^{d}(\mathbb{P}^n)$. 

Let us first determine the vertices of the $A$-graph of $Hilb^{d}(\mathbb{P}^n)$. Any vertex must be $T$-fixed, therefore it must be a monomial ideal. The $G_{a}$-fixed monomial ideals agree with the Borel fixed ideals (\cite{Eis}, thm. 15.23): 

\begin{lem} \label{lem:vertices}
A homogenous monomial ideal $I$ is $G_{a}$-fixed if and only if it has the following property: For every $f\in I$ divisible by $X_{i}$ one has $\frac{X_{i+1}}{X_{i}}f\in I$. 
\end{lem} 

\noindent \textit{Proof} Suppose that $I$ satisfies the given condition. Since $\theta(a)(X_{i})=X_{i}+aX_{i+1}+\ldots+\frac{a^{n-i}}{(n-i)!}X_{n}$ it is clear that if $m$ is a monomial in $I$ then every monomial summand of $\theta(a)m$ remains in $I$. Therefore $I$ is $G_{a}$-fixed. Conversely, suppose that $I$ is $G_{a}$-fixed. Say $f\in I$ and $X_{i}|f$. Without loss of generality, we may assume that $f=\prod X_{j}^{e_{j}}$ is a monomial. Then $\theta(a)(f)=\prod (X_{j}+aX_{j+1}+\ldots)^{e_{j}}\in I$. Since $e_{i}\geq 1$, by binomial expansion the monomial $\frac{X_{i+1}}{X_{i}}f$ appears among the monomial summands of $\theta(a)(f)$. Since $I$ is a monomial ideal, $\frac{X_{i+1}}{X_{i}}f\in I$. $\Box$    

Let us now pass to the affine chart $X_{0}\neq 0$, which contains the fixed point $[1:0:\ldots:0]$. Say $\displaystyle{x_{i}=\frac{X_{i}}{X_{0}}}$ be affine coordinates on this chart. Since the support of a $G_{a}$-fixed monomial ideal in the Hilbert scheme must be the point $[1:0:\ldots:0]$, without loss of generality we may consider monomial ideals in the variables $x_{1},\ldots,x_{n}$ from now on.  

\begin{defn}
Suppose that $I$ and $J$ are two $G_{a}$-fixed monomial ideals. Say $m_{1},m_{2},\ldots,m_{k}$ and $m_{1}^{\prime},\ldots,m_{k}^{\prime}$ are monomials such that $m_{i}<m_{i}^{\prime}$ in lexicographic order for each $i$. Assume that each $m_{i}$ and $m_{i}^{\prime}/x_{1}$ belongs to the minimal set of generators of $I$. We say that $(m_{1},\ldots,m_{k})\rightarrow (m_{1}^{\prime},\ldots,m_{k}^{\prime})$ is a move from $I$ to $J$ if $J$ can be obtained from $I$ by replacing $m_{i}$ by $m_{i}/x_{1}$ and $m_{i}^{\prime}/x_{1}$ by $m_{i}^{\prime}x_{j}/x_{1}$ for $j=1,\ldots,n$ in the minimal set of 
generators of $I$. 
\end{defn}

\begin{rem} After a move, the set of generators obtained for $J$ need not be minimal anymore. 
\end{rem}

\begin{rem} \label{rem:socle}
A move from $I$ to $J$ can easily be expressed in terms of socles of these ideals: The move $(m_{1},\ldots,m_{k})\rightarrow (m_{1}^{\prime},\ldots,m_{k}^{\prime})$ removes $m_{i}/x_{1}$ from the socle of $I$ and adds $m^{\prime}_{i}/x_{1}$ to it for each $i=1,\ldots,k$. Such an operation does not change the number of elements in the socle. Therefore, if there exists a move from $I$ to $J$ then $I$ and $J$ must have the same Hilbert polynomial. 
\end{rem} 

\begin{ex} Say $n=3$. Suppose that 
\begin{align*}I&=<x_{1}^{3},x_{1}^{2}x_{2}, x_{1}^{2}x_{3}, x_{1}x_{2}^{2}, x_{1}x_{2}x_{3}, x_{1}x_{3}^{2}, x_{2}^{3}, x_{2}^{2}x_{3}, x_{2}x_{3}^{2}, x_{3}^{3}>\\ 
J&=<x_{1}^{4}, x_{1}^{2}x_{2}, x_{1}x_{2}^{2}, x_{1}^{2}x_{3}, x_{1}x_{2}x_{3}, x_{2}^{3}, x_{2}^{2}x_{3}, x_{3}^{2}>.
\end{align*} 
Take $m_{1}=x_{1}x_{3}^{2}$ and $m_{1}^{\prime}=x_{1}^{4}$. Then $J$ can be obtained from $I$ by the move $m_{1}\rightarrow m_{1}^{\prime}$. 
\end{ex}

\begin{defn}
Say $m=x_{1}^{a_1}x_{2}^{a_2}\ldots x_{n}^{a_{n}}$ is a monomial. Define the weight of $m$ to be 
\[ w(m)=\sum_{j=1}^{n} a_{j} (n-j) \]
The weight of a monomial ideal $I$ is defined to be the sum of the weights of all monomials in the socle of $I$. 
\end{defn}

\begin{lem} \label{lem:weight} Suppose that $I$ and $J$ are two $G_{a}$-fixed monomial ideals. If there exists a move from $I$ to $J$ then $w(J)>w(I)$. 
\end{lem}

\noindent\textit{Proof} By remark \ref{rem:socle}, it suffices to show that $w(m_{i}/x_{1})<w(m_{i}^{\prime}/x_{1})$, namely $w(m_{i})<w(m_{i}^{\prime})$. Since $m_{i}<m_{i}^{\prime}$ in lexicographic ordering and $I$ is $G_{a}$-fixed, $deg(m_{i})\leq deg(m_{i}^{\prime}/x_{n})<deg(m_{i}^{\prime})$. The result follows. $\Box$

Let $I$ be a $G_{a}$-fixed monomial ideal minimally generated by $S=\{m_{1},m_{2},\ldots,m_{k}\}$. Let $S^{(d)}$ denote the subset of degree $d$ monomials in $S$. 

\begin{defn}
Let $d$ be the largest integer such that $S^{(d)}\neq \emptyset$ and $S^{(d)}$ contains a monomial other than $x_{1}^{d}$. Let $j$ be the maximal positive integer such that $x_{j}$ divides some element of $S^{(d)}$. Let $l$ be the maximal integer such that $x_{j}^{l}$ divides at least two elements of $S^{(d)}$. Let 
\[ \hat{S}=\{m_{i}\in S^{(d)}: x_{j}^{l}|m_{i}, x_{1}|m_{i}  \} \]
\end{defn}

Let us write elements of $\hat{S}$ in increasing order with respect to lex. Namely, let $\hat{S}=\{m_{1},m_{2},\ldots,m_{r}\}$ such that $m_{1}<m_{2}<\ldots<m_{r}$. 

\begin{lem} \label{lem:multiple}
Suppose that $I$, $j$, $l$ and $\hat{S}=\{m_{1},\ldots,m_{r}\}$ are as above and that $r\geq 2$. Then there exists a positive integer $s$ such that the ideal $J$ obtained from $I$ by the move 
\[ (m_{1}, m_{2}, \ldots, m_{s})\mapsto (\left(\frac{x_{1}}{x_{j}}\right)^{l}x_{1}m_{1},\left(\frac{x_{1}}{x_{j}}\right)^{l}x_{1}m_{2},\ldots,\left(\frac{x_{1}}{x_{j}}\right)^{l}x_{1}m_{s}) \]
is also $G_{a}$-fixed. 
\end{lem}

\noindent \textit{Proof} Since the largest element in $\hat{S}$ is $x_{1}^{d-l}x_{j}^{l}$ and $r\geq 2$, we have $m_{1}\neq x_{1}^{d-l}x_{j}^{l}$. Then $m_{1}$ must be divisible by $x_{i}$ for some $1<i<j$. By definition of $\hat{S}$, the monomial $\left(\frac{x_{1}}{x_{j}}\right)^{l}m_{1}$ must be in $I$, hence the move $m_{1}\mapsto x_{1}\left(\frac{x_{1}}{x_{j}}\right)^{l}m_{1}$ is valid. If the resulting ideal is not $G_{a}$-fixed then $\frac{x_{1}}{x_{i}}\left(\frac{x_{1}}{x_{j}}\right)^{l}m_{1}$ must be in $I$ for some $1<i<j$. Select the maximal value of $i$ for which this holds. Then by $G_{a}$ invariance of $I$, the monomial $\frac{x_{1}}{x_{i}}m_{1}=m_{k}$ must be in $\hat{S}$. Apply the move 
\[(m_{1}, m_{2}, \ldots, m_{k})\mapsto (\left(\frac{x_{1}}{x_{j}}\right)^{l}x_{1}m_{1},\left(\frac{x_{1}}{x_{j}}\right)^{l}x_{1}m_{2},\ldots,\left(\frac{x_{1}}{x_{j}}\right)^{l}x_{1}m_{k}). \]
Repeat the same argument with $m_{1}$ replaced by $m_{k}$. After finitely many steps, the process terminates. $\Box$ 

\begin{lem} \label{lem:single} 
Suppose that $I,j,l$ and $\hat{S}$ are as above but $r=1$, so that $\hat{S}=\{m_1\}$. Then there exists $h\leq l$ and $k\geq h$ such that the ideal $J$ obtained from $I$ by the move 
\[ m_{1}\mapsto \left( \frac{x_{1}^{k}}{x_{j}^{h}} \right) x_{1}m_{1} \]
is also $G_{a}$-fixed. If $h<l$ then $k=h$. 
\end{lem}
 
\noindent \textit{Proof} First suppose that $l\geq 2$. If the monomial $\displaystyle{ \frac{x_{1}^{l}}{x_{j}^{l}}m_{1}}$ is in $I$ then the move  
\[ m_{1}\mapsto \left( \frac{x_{1}^{l}}{x_{j}^{l}} \right) x_{1}m_{1} \]
is valid and the resulting ideal is $G_{a}$-fixed. If not, then by definition of $\hat{S}$, the monomial $\displaystyle{ \frac{x_{1}^{l-1}}{x_{j}^{l-1}}m_{1}}$ is in $I$ so the move 
\[ m_{1}\mapsto \left( \frac{x_{1}^{l-1}}{x_{j}^{l-1}} \right) x_{1}m_{1} \]
is valid and again the resulting ideal is $G_{a}$-fixed. 

It remains to consider the case $l=1$. Then $m_{1}=x_{1}^{d-l}x_{j}^{l}$. Suppose that $d^{\prime}$ is minimal such that $x_{1}^{d^{\prime}}\in I$. Then the move 
\[m_{1}\mapsto \left( \frac{x_{1}^{d^{\prime}-d+l}}{x_{j}^{l}} \right) x_{1}m_{1} \]
is valid and again, the ideal obtained is $G_{a}$-fixed. 
$\Box$

Lemmas \ref{lem:weight}, \ref{lem:multiple} and \ref{lem:single} imply that starting from any $G_{a}$-fixed monomial ideal $I$ there exists a finite sequence of moves through the vertices of the $A$-graph of the Hilbert scheme ending at the ideal $<x_{1}^{d},x_{2},\ldots,x_{n-1},x_{n}>$. Next we want to show that if we can pass from vertex $v$ to vertex $w$ by a move as described above, then these two vertices are connected with an edge of the $A$-graph. 

\begin{lem} \label{lem:edges} Suppose that $I$ is minimally generated by $m_1,m_2,\ldots, m_{k}$ and $\hat{S}=\{m_{1},\ldots,m_{r}\}$, $s$, $J$, $j$ and $l$ are as in lemma \ref{lem:multiple}. Assume that $J$ is obtained from $I$ by the move 
\[ (m_{1},m_{2},\ldots,m_{s})\mapsto (m_{1}^{\prime}, m_{2}^{\prime}, \ldots,m_{s}^{\prime}). \]
Reorder the minimal generating set for $I$ such that the first $s$ elements are $m_{1}^{\prime}/x_{1},\ldots,m_{s}^{\prime}/x_{1}$ and the remaining ones are $\hat{m}_{s+1},\ldots, \hat{m}_{k}$. Suppose that $x_{1}^{a_{i}}| m_{i}$ but $x_{1}^{a_{i}+1}$ does not divide $m_{i}$.

Then,

(i) The family of ideals 
\[ I_{t}= \left< \frac{m_{1}^{\prime}}{x_{1}}+c_{1}t\frac{m_{1}}{x_{1}}, \frac{m_{2}^{\prime}}{x_{1}}+c_{2}t\frac{m_{2}}{x_{1}},\ldots,\frac{m_{s}^{\prime}}{x_{1}}+c_{s}t\frac{m_{s}}{x_{1}}, \hat{m}_{s+1},\ldots,\hat{m}_{k} \right> \]
where $\displaystyle{c_{i}= \frac{(a_{i}+l)!}{(a_{i}-1)!} }$ forms an edge of the $A$-graph.

(ii) $I_{0}=I$. 

(iii) $\lim_{t\rightarrow \infty} I_{t}=J$. 
\end{lem}

\noindent\textit{Proof} Assuming that $I_{t}$ is defined as in (i), the proof of (ii) immediately follows by putting $t=0$. 

In order to prove (iii) it suffices to show that $m_{i}^{\prime}\in J$ and $\frac{m_{i}}{x_{1}}\in J$ for each $i\in \{1,\ldots,s\}$. 
Since $m_{i}\in I_{t}$ for each $t$, the first of these follows from the equality 
\[ m_{i}^{\prime}=x_{1}\left( \frac{m_{i}^{\prime}}{x_{1}}+c_{i}t\frac{m_{i}}{x_{1}}\right)-c_{i}tm_{i}. \] 
The second can be shown by dividing $\frac{m_{i}^{\prime}}{x_{1}}+c_{i}t\frac{m_{i}}{x_{1}}$ by $c_{i}t$ and taking $t$ to $\infty$. 

Now let us prove (i). We need to show that $I_{t}$ is $G_{a}$-fixed for each  $t$. 
Let us write the generators $\{\hat{m}_{s+1},\ldots,\hat{m}_{k}\}$ in decreasing order. The last monomial, which is of the form $x_{n}^p$, is $G_{a}$-fixed. The $G_{a}$-action on a monomial $m$ produces monomials less than or equal to $m$ in the lexicographic order. Using this fact and descending induction on the index, we see that $<\hat{m}_{s+1},\ldots,\hat{m}_{k}>$ is $G_{a}$-fixed. Similarly we see that $< \frac{m_{1}^{\prime}}{x_{1}}+c_{1}t\frac{m_{1}}{x_{1}}, \hat{m}_{s+1},\ldots,\hat{m}_{k}>$ is $G_{a}$-fixed. We want to show by induction that each segment $< \frac{m_{1}^{\prime}}{x_{1}}+c_{1}t\frac{m_{1}}{x_{1}}, \frac{m_{2}^{\prime}}{x_{1}}+c_{2}t\frac{m_{2}}{x_{1}},\ldots,\frac{m_{i}^{\prime}}{x_{1}}+c_{i}t\frac{m_{i}}{x_{1}}, \hat{m}_{s+1},\ldots,\hat{m}_{k}>$ is $G_{a}$-fixed. In order to complete the induction step, we need to see that the $G_{a}$-action on $\frac{m_{i}^{\prime}}{x_{1}}+c_{i}t\frac{m_{i}}{x_{1}}$ yields the monomials $\frac{m_{r}^{\prime}}{x_{r}}$ and $\frac{m_{r}}{x_{r}}$ in a term proportional to $\frac{m_{r}^{\prime}}{x_{1}}+c_{r}t\frac{m_{r}}{x_{1}}$. First assume $a_{i}=a_{r}$, so that $c_{i}=c_{r}$. Then the coefficient of $\frac{m_{r}}{x_{1}}$ resulting from the $G_{a}$-action on $\frac{m_{i}}{x_{1}}$ is a weighted sum over all decreasing sequences (in lexicographic order) of monomials starting from $\frac{m_{i}}{x_{1}}$ and ending at $\frac{m_{r}}{x_{1}}$. Since the move $m_{i}\rightarrow m_{i}^{\prime}$ does not change any exponents other than $x_{1}$ and $x_{j}$, the same set of decreasing sequences occurs from $\frac{m_{i}^{\prime}}{x_{1}}$ to $\frac{m_{r}^{\prime}}{x_{r}}$. Therefore the aformentioned claim holds. On the other hand if $a_{i}\neq a_{r}$, the only difference in the sequences is in the steps transferring powers of $x_{1}$ to the other variables, and the choice of $c_{i}$ balances the count. 
$\Box$ 

\begin{cor} The ideals in lemma \ref{lem:vertices} correspond to the vertices and the families of ideals in lemma \ref{lem:edges} correspond to edges of a spanning subtree of the $A$-graph of $Hilb^{d}(\mathbb{P}^{n})$. In particular, this graph is connected, hence by theorem \ref{thm:Agraph}, $Hilb^{d}(\mathbb{P}^{n})$ is connected. 
\end{cor}

\end{document}